\documentclass[12pt,reqno]{amsproc}

\newtheorem{theo}{Theorem}

\newtheorem{lem}{Lemma}

\newtheorem{prop}{Proposition}

\newcommand\eps\varepsilon
\newcommand\ph\varphi
\newcommand\kap\Lambda

\usepackage{enumerate}
\usepackage{graphicx}
\usepackage{float}

\begin{document}

\title[Br\o ndsted's Fixed Point Theorem]
{A Note on Br\o ndsted's Fixed Point Theorem}

\author[Oleg Zubelevich]{Oleg Zubelevich\\
 Steklov Mathematical Institute of Russian Academy of Sciences \\oezubel@gmail.com
 }

\date{}
\thanks{The research was funded by a grant from the Russian Science
Foundation (Project No. 19-71-30012)}

\begin{abstract}We show that for the case of uniformly convex Banach spaces
the conditions of the Br\o ndsted fixed point theorem can be relaxed.
\end{abstract}

\maketitle
\numberwithin{equation}{section}
\newtheorem{theorem}{Theorem}[section]
\newtheorem{lemma}[theorem]{Lemma}
\newtheorem{definition}{Definition}[section]

\section{Introduction. The main theorem}
The object of this short note is a fixed point theorem by Arne Br\o ndsted.
Let us formulate this theorem.

Let $(X,\|\cdot\|)$ be a Banach space and let $M\subset X$ be a closed set.
We denote a closed unit ball as $B=\{\|x\|\le 1\}$.
Assume that
\begin{equation}\label{asdf44}
M\cap B=\emptyset.\end{equation}
Consider a mapping $T:M\to M$ that maps each $x\in M$ in the direction of the ball:
if $Tx\ne x$ then there exists $t>1$ such that
\begin{equation}\label{sdfsgtt}x+t(Tx-x)\in B.\end{equation}
\begin{theo}[\cite{arne}]\label{dsfg545} In addition to the assumptions above suppose that
\begin{equation}\label{sfrg5500}
\inf\{\|x\|\mid x\in M\}>1.\end{equation}Then the mapping $T$ has a fixed point.\end{theo}
Observe that condition (\ref{sfrg5500}) is stronger than condition (\ref{asdf44}) if only $\dim X=\infty$.

To prove theorem \ref{dsfg545} Br\o ndsted endows the set $M$ with a partial order in the following way.

{\it If $x,y\in M$ then by definition we write $x<y$ provided either $x=y$ or there exists $t>1$ such that
$$x+t(y-x)\in B.$$} Formula (\ref{sdfsgtt}) takes the form
\begin{equation}\label{sdf44}x<Tx,\quad \forall x\in M.\end{equation}
Then Br\o ndsted observes that this partial order is finer than one of the Caristi type \cite{car} and   thus by some other of his results \cite{arbe2} the fixed point exists.

Our aim is to show that for the class of uniformly convex Banach spaces $X$ theorem \ref{dsfg545} remains valid even in the critical case  when  condition (\ref{sfrg5500}) is replaced with  (\ref{asdf44}). This fact does not follow from the original Br\o ndsted's method.

Recall a definition. {\it A Banach space $(X,\|\cdot\|)$ is said to be uniformly convex if for any $\sigma>0$ there exists $\gamma>0$ such that if $\|x\|=\|y\|=1$ and $\|x-y\|\ge \sigma$ then $\|x+y\|\le 2-\gamma.$}

For example the spaces $L^p,\quad p\in(1,\infty)$ are uniformly convex; each uniformly convex Banach space is reflexive; a Hilbert space is uniformly convex, see \cite{yos} and references therein.
\begin{theo}\label{sxfg60oo}Assume that $X$ is a uniformly convex Banach space.
If the mapping $T$ satisfies condition (\ref{sdf44}) and condition (\ref{asdf44}) is fulfilled then $T$ has a fixed point.\end{theo}
\section{Proof of theorem \ref{sxfg60oo}}
The  scheme of the proof is quite standard by itself.  It is clear that a maximal element of the set $M$ provides a fixed point. To prove that the maximal element exists we check the conditions of the Zorn lemma.

This argument and the technique developed below give possibility to conduct a direct proof of theorem \ref{dsfg545} as well.

\begin{prop}\label{dsfg55}Suppose that vectors $a,x\in X$ have the following properties $$\|(1-t)a+tx\|>1,\quad \forall t\in (0,1),\quad \|x\|>1,\quad \|a\|=1.$$ Then for any $\eps>0$ there exists $\delta>0$ such that
inequality $\|x\|\le 1+\delta$ implies $\|x-a\|\le \eps.$ \end{prop}
This proposition admits a ''physical'' interpretation. Let $x$ be a light source placed away from the ball $B:\quad \|x\|>1$. According to the proposition the diameter of the light spot on the ball tends to zero as $x$ approaches  the ball: $\|x\|\to 1$.

Here the uniform convexity of the norm is essential:
such a feature is not yet kept  for the norm $\|(p,q)\|=\max\{|p|,|q|\}$ in $\mathbb{R}^2$.

{\it Proof of proposition \ref{dsfg55}.}
Assume the opposite: there exist  $\eps>0$ and  sequences $a_n,x_n,$ \begin{equation}\label{gkk}\|x_n\|>1,\quad\|a_n\|=1,\quad  \|x_n\|\to 1,\quad\|(1-t)a_n+tx_n\|>1 \end{equation} such that
$$\|x_n-a_n\|> \eps.$$
Consequently, for all sufficiently large $n$  the estimate
\begin{equation}\begin{split}
\|x_n&-a_n\|=\Big\|x_n-\frac{x_n}{\|x_n\|}+\frac{x_n}{\|x_n\|}-a_n\Big\|\nonumber\\
&\le\alpha_n+\Big\|a_n-\frac{x_n}{\|x_n\|}\Big\|,\quad \alpha_n=\|x_n\|\Big(1-\frac{1}{\|x_n\|}\Big)\to 0\end{split}\end{equation}
implies
$$\Big\|a_n-\frac{x_n}{\|x_n\|}\Big\|\ge\eps/2.$$

Plugging  $t=1/2$ in (\ref{gkk}) we obtain
\begin{equation}\label{sdfg0hhjuy}\|a_n+x_n\|>2.\end{equation}
The inequality
$$\Big\|a_n+\frac{x_n}{\|x_n\|}\Big\|>2-\alpha_n$$ follows from (\ref{sdfg0hhjuy}) in the same way as above.

This contradicts  the hypothesis of uniform convexity of the space $X$.

The proposition is proved.

Let $C\subset M$ be a chain then put $\rho=\inf\{\|u\|\mid u\in C\},\quad \rho\ge 1$.

The inclusion $x\in C$ implies that $\|x\|>1$ provided $\rho=1$ and $\|x\|\ge \rho$ provided $\rho>1$.

For any $x\in C$ define a set
$$K_x(\rho)=\{y\in M\mid \|y\|\ge \rho,\quad y>x\}.$$
The sets $K_x(\rho)$ are  nonvoid: $x\in K_x(\rho)$ and
\begin{equation}\label{0}x_1<x_2\Longrightarrow K_{x_2}(\rho)\subset K_{x_1}(\rho).\end{equation}
\begin{lem}\label{sxdfgoop}
 The sets $K_x(\rho)$ are closed. \end{lem}{\it Proof of lemma \ref{sxdfgoop}.} Indeed, let a convergent sequence $\{y_k\}$ belong to $K_x(\rho)$ and $y_k\to y\in M$.
This means that there are  sequences $\{\beta_k\}\subset(0,1)$ and  $\{a_k\}\subset X,\quad \|a_k\|=1$ such that
$$y_k=\beta_k a_k+(1-\beta_k)x.$$ The sequence $\{\beta_k\}$ contains a convergent subsequence; we keep the same notation for this subsequence: $\beta_k\to \beta$.

If $\beta=0$ then $\|\beta_ka_k\|\to 0$ and $y=x\in K_\rho(x)$. If $\beta\ne 0$ then
$$a_k=\frac{1}{\beta_k}y_k+\Big(1-\frac{1}{\beta_k}\Big)x\to a,\quad \|a\|=1$$ and
$$y=\beta a+(1-\beta)x.$$
Observe that since $y\in M$ the parameter $\beta$ cannot be equal to $1$.

The lemma is proved.

\begin{lem}\label{2}Suppose that $z\in K_x(\rho),\quad x\in C$.

If $\rho> 1$ then
$$\|z-x\|\le(\|x\|-\rho)\frac{\rho+1}{\rho-1}.$$
If $\rho=1$ then for any $\eps>0$ there exists  $\delta>0$ such that
$$\|x\|\le 1+\delta\Longrightarrow \|z-x\|\le\eps.$$
\end{lem}
{\it Proof of lemma \ref{2}.}

{\bf The case $\rho>1$.}
Indeed, the formulas \begin{equation}\label{dxfh600}
x+t(z-x)=a,\quad \|a\|=1,\quad t>1,\quad \|x\|,\|z\|\ge \rho>1\end{equation} imply
$z=(a+(t-1)x)/t$ and
$$\rho\le \|z\|\le \frac{1}{t}+\frac{t-1}{t}\|x\|,\quad \frac{1}{t}\le \frac{\|x\|-\rho}{\|x\|-1}.$$
Use (\ref{dxfh600}) again:
$$\|z-x\|= \frac{1}{t}\|a-x\|\le\frac{1}{t}(1+\|x\|).$$
{\bf The case $\rho=1$.}
The condition of the lemma $z\in K_x(1)$ means
$$z=\tau a+(1-\tau) x,\quad \|x\|,\|z\|>1,\quad \|a\|=1,\quad \tau\in(0,1).$$
Therefore the assertion of the lemma follows from proposition \ref{dsfg55} and the formulas
$$z-x=\tau(a-x),\quad \|z-x\|<\|x-a\|.$$

The lemma is proved.

\begin{lem}\label{3} For any $\eps>0$ there exists $\tilde x\in C$ such that
$$C\ni x>\tilde x\Longrightarrow \mathrm{diam}\,K_x(\rho)<\eps.$$\end{lem}
{\it Proof of lemma \ref{3}.}
 Consider the case $\rho>1$. By definition of the number $\rho$, for any $\eps>0$ there is an element $\tilde x\in C$ such that
$$ \|\tilde x\|\le\eps+\rho.$$
Take any elements $z_1,z_2\in K_{\tilde x}$ and apply lemma \ref{2} for each summand from the right side  of the inequality
$$\|z_1-z_2\|\le \|z_1-\tilde x\|+\|z_2-\tilde x\|.$$
Observe also that  formula (\ref{0}) implies $$\tilde x<x\in C\Longrightarrow \mathrm{diam}\,K_{x}(\rho)\le \mathrm{diam}\,K_{\tilde x}(\rho).$$

The case $\rho=1$ is processed analogously be means of lemma \ref{2}.

The lemma is proved.

Therefore we have a nested family of the closed sets $K_x(\rho)$ which diameters tend to zero.
By well-known theorem their intersection is not empty and consists of a single point:
$$\bigcap_{x\in C}K_x(\rho)=\{m\}.$$ The point  $m\in M$ is an upper bound for $C$. Indeed, for any $x\in C$ we have $m\in K_x(\rho)$ and thus $x<m$.

The theorem is proved.

\end{document}